\documentclass[preprint,12pt]{elsarticle}

\usepackage{amssymb}
\usepackage{amsmath}
\usepackage{amsthm}

\newtheorem{theorem}{Theorem}%
\newtheorem{proposition}[theorem]{Proposition}%

\newtheorem{lemma}[theorem]{Lemma}

\newtheorem*{definition}{Definition}

\newtheorem*{conjecture}{Conjecture}

\journal{Finite Fields and Their Applications}

\begin{document}

\begin{frontmatter}



\title{On an Erd\H{o}s-type conjecture on $\mathbb{F}_q[x]$}


\author{Rongyin Wang} 

\affiliation{organization={School of Mathematical Sciences and LPMC, Nankai University},
            addressline={94 Weijin Road}, 
            city={Tianjin},
            postcode={300071}, 
            country={China}}

\begin{abstract}
P. Erd\H{o}s conjectured in 1962 that on the ring $\mathbb{Z}$, every set of $n$ congruence classes in $\mathbb{Z}$ that covers the first $2^n$ positive integers also covers the ring $\mathbb{Z}$. 
This conjecture was first confirmed in 1970 by R. B. Crittenden and C. L. Vanden Eynden. 
Later, in 2019, P. Balister, B. Bollob\'{a}s, R. Morris, J. Sahasrabudhe, and M. Tiba provided a more transparent proof. 
In this paper, we follow the approach used by R. B. Crittenden and C. L. Vanden Eynden to prove the generalized Erd\H{o}s' conjecture in the setting of polynomial rings over finite fields. 
We prove that every set of $n$ cosets of ideals in $\mathbb F_q[x]$ that covers all polynomials whose degree is less than $n$ covers the ring $\mathbb{F}_q[x]$. 
\end{abstract}


\end{frontmatter}



\section{Introduction}\label{introsec}

In 1970, R. B. Crittenden and C. L. Vanden Eynden\cite{CV1970} proved that any set of $n$ congruence classes covering the positive integers $1, 2, \cdots, 2^n$ covers all integers, which was originally conjectured by P. Erd\H{o}s\cite{Erdos1962} in 1962, who first proved a weaker version of this conjecture, replacing $2^n$ with $n \cdot 2^n$. It is easy to see that the constant $2^n$ cannot be improved, since one can find $n$ congruence classes that cover any given $2^n$ consecutive integers with any one given exception among them exactly. In 2019, P. Balister, B. Bollob\'as, R. Morris, J. Sahasrabudhe, and M. Tiba\cite{BBST2020} gave a simpler proof of this conjecture by Erd\H{o}s, including the construction of polynomials using parameters derived from counterexample congruence classes.

Covering systems for polynomial rings over finite fields have been studied by \cite{Azlin2011thesis, LWWY2024, LWWY2024_2}. The main result of this paper is that we prove the generalized Erd\H{o}s' conjecture in the ring $\mathbb{F}_q[x]$. We define covering systems in polynomial rings in a similar manner to how they are defined in the integer ring. 
\begin{definition}
In a polynomial ring $\mathbb{F}_q[x]$ over a finite field $\mathbb{F}_q$, a congruence class, denoted by $\langle f(x) \rangle + r(x)$ or $r(x) \pmod{f(x)}$, is the following set of polynomials
\[
\{
g(x) \in \mathbb{F}_q[x] : g(x) \equiv r(x) \pmod{f(x)}
\}. 
\]
We say that a set of congruence classes covers a certain subset $A \subseteq \mathbb{F}_q[x]$ if, for any $g(x) \in A$, $g(x)$ belongs to at least one of the congruence classes. In the case $A = \mathbb{F}_q[x]$, we say that the set of congruence classes is a covering system of $\mathbb{F}_q[x]$. 
\end{definition} 
The main result of this paper is the following theorem. 
\begin{theorem}\label{2thm}
Let $f_1(x), f_2(x), \cdots, f_n(x), r_1(x), r_2(x), \cdots, r_n(x)$ be polynomials in $\mathbb{F}_2[x]$ satisfying $\deg(r_i(x)) < \deg(f_i(x))$ for $1\leq i\leq n$. Suppose that there exists some polynomial $g(x) \in \mathbb{F}_2[x]$ that satisfies none of the congruence equations
\[
g(x) \equiv r_i(x) \pmod{f_i(x)}, 
\]
for all $1 \leq i \leq n$, then there exists a polynomial $g_0(x) \in \mathbb{F}_2[x]$ that satisfies none of the congruence equations above and $\deg(g_0(x)) < n$. 
\end{theorem}
In fact,  our result is also valid for $\mathbb{F}_q[x]$ employing the same proof methodology. 
\begin{theorem}\label{qthm}
Let $f_1(x), f_2(x), \cdots, f_n(x), r_1(x), r_2(x), \cdots, r_n(x)$ be polynomials in $\mathbb{F}_q[x]$ satisfying $\deg(r_i(x)) < \deg(f_i(x))$ for $1\leq i\leq n$. Suppose that there exists some polynomial $g(x) \in \mathbb{F}_q[x]$ that satisfies none of the congruence equations
\[
g(x) \equiv r_i(x) \pmod{f_i(x)}, 
\]
for all $1 \leq i \leq n$, then there exists a polynomial $g_0(x) \in \mathbb{F}_q[x]$ which satisfies none of the congruence equations above and $\deg(g_0(x)) < n$. 
\end{theorem}
We believe that the result in Theorem~\ref{qthm} is not sharp. We conjecture that one can replace the degree of the polynomial $\deg(g_0(x)) < n$ with $\deg(g_0(x)) < {n}/{(q - 1)}$. 
\begin{conjecture}
Let $f_1(x), f_2(x), \cdots, f_n(x), r_1(x), r_2(x), \cdots, r_n(x)$ be given polynomials in $\mathbb{F}_q[x]$ satisfying $\deg(r_i(x)) < \deg(f_i(x))$ for $1\leq i\leq n$. Suppose that there exists some polynomial $g(x) \in \mathbb{F}_q[x]$ that satisfies none of the congruence equations
\[
g(x) \equiv r_i(x) \pmod{f_i(x)}, 
\]
for all $1 \leq i \leq n$, then there exists a polynomial $g_0(x) \in \mathbb{F}_q[x]$ which satisfies none of the congruence equations above and $\deg(g_0(x)) < {n}/{(q - 1)}$. 
\end{conjecture}
Our proof mainly adheres to the method established by R. B. Crittenden and C. L. Vanden Eynden \cite{CV1970}. The reason we do not adopt the more straightforward proof by P. Balister et al. \cite{BBST2020} lies in two critical aspects: first, their approach hinges on mapping the elements of the ring to an algebraically closed field while preserving divisibility, and second, it involves identifying specific monomials within this field that exhibit linear independence. Extending their method to polynomial rings over finite fields poses challenges in achieving these objectives.

The paper is structured as follows. 
In Section~\ref{keypropsec}, we lay out some essential propositions related to our theorems, showing that our theorem is applicable for sufficiently large values of $n$ and that our result for $\mathbb{F}_2[x]$ is indeed sharp. 
In Section~\ref{lemmassec}, we articulate and establish the essential lemmas to prove Proposition~\ref{keyprop}. 
Section~\ref{proofsec} is dedicated to the proof of Proposition~\ref{keyprop}, employing the adjustment method. In Section~\ref{futureresultssec}, we investigate the scenario of Theorem~\ref{2thm} for the cases where $n \leq 9$. Finally, in Section~\ref{generalizedresultssec}, we extend our findings to the broader context of $\mathbb{F}_q[x]$. 

\section{Observations for $\mathbb{F}_2[x]$}~\label{keypropsec}
In this section, we discuss some observations we have for the polynomial ring $\mathbb{F}_2[x]$.

The first proposition we would introduce is one that played a key role in our proof of Theorem\ref{2thm}.
\begin{proposition}\label{keyprop}
If Theorem~\ref{2thm} is false, it does not hold for some $n<10$. In other words, there exists $n<10$ that we can find $n$ congruence equations covering all polynomials $g(x)\in \mathbb{F}_2[x]$ with $\deg(g(x))<n$ but not the whole ring.
\end{proposition}

We will prove this proposition in the following sections, in particular Sections \ref{lemmassec} and \ref{proofsec}. 

Another proposition of the ring shows the sharpness of our Theorem \ref{2thm}, i.e., the degree $n$ in Theorem \ref{2thm} cannot be replaced by $n-1$. Actually, these polynomials acts like the $2^n$ consecutive integers in the original problem, and sets of congruences which fail to cover the whole ring may omit any given one among them.

The idea here is that, just like how one can find a set of $n$ congruence equations covering all but any given integer in $[1,2^n]$ (using the bisection method), we find a set of $n$ congruence equations covering all but any given polynomial whose degree is at most $n-1$.

\begin{proposition}
    Let $n\in\mathbb{N}$, $r_0(x) \in \mathbb{F}_2[x]$ be any given polynomial with $\deg(r_0(x))<n$, one can find a set of $n$ congruence equations covering all polynomials $g(x)$ with $\deg(g(x))<n$, except for $r_0(x)$.
\end{proposition}
\begin{proof}
 We start by finding the residues of $r_0(x)$ modulo $x,x^2,\cdots x^n$ respectively. So there exist $r_1(x),r_2(x), \cdots,r_n(x)$ such that 
\[
\begin{cases}
r_0(x) \equiv r_1(x) \pmod{x}, \\
r_0(x) \equiv r_2(x) \pmod{x^2}, \\
\cdots \\
r_0(x) \equiv r_n(x) \pmod{x^n}. 
\end{cases}
\]
One can find $n$ congruence classes, starting with $r'_1(x) \equiv 1+r_1(x) \pmod{x}$. There are 4 elements in the quotient ring $\mathbb{F}_2[x]/\langle x^2 \rangle$ and $r_2(x) \not\equiv 1+r_1(x) \pmod{x}$, so there is exactly one $r'_2(x)$ whose degree is less than 2, $r'_2(x) \not\equiv r_2(x) \pmod{x^2}$ and $r'_2(x) \not\equiv 1+r_1(x) \pmod{x}$.

Suppose we have already in this way found $r'_1(x) (=1+r_1(x)), r'_2(x), \cdots, r'_k(x)$. By the way they are constructed we know 
$\langle x \rangle +r'_1(x), 
\langle x^2 \rangle +r'_2(x), 
\cdots, 
\langle x^k \rangle +r'_k(x), 
\langle x^{k+1} \rangle + r_{k+1}(x)$
are disjoint in $\mathbb{F}_2[x]$, hence 
${\langle x \rangle}+{r'_1(x)}$, 
${ \langle x^2 \rangle}+{r'_2(x)}$, 
$\cdots$, 
${\langle x^k \rangle}+{r'_k (x)}$ and $\{{r_{k+1}(x)}\}$ are disjoint in $\mathbb{F}_2[x]/\langle x^{k+1} \rangle$. 
Since there are $2^{k+1-d}$ elements of $\mathbb{F}_2[x]/\langle x^{k+1} \rangle$ in ${\langle x^d \rangle}+{r'_d (x)}$, we know there is exactly one $r'_{k+1}(x) \in \mathbb{F}_2[x]$, such that $\deg(r'_{k+1}(x)) < k+1$, that $\langle x \rangle+r'_1(x), 
\langle x^2 \rangle+r'_2(x),
\cdots, 
\langle x^k \rangle +r'_k(x), 
\langle x^{k+1} \rangle +r_{k+1}(x),
\langle x^{k+1} \rangle +r'_{k+1}(x)$ disjoint, $\left(
\bigcup\limits_{i=1}^{k+1}
( \langle x^i \rangle +r'_i(x) ) \cup (\langle x^{k+1} \rangle + r_{k+1}(x))
\right)
= \mathbb{F}_2[x]$. In other words, the set of congruence classes:
\[
\{
r'_1(x) \pmod{x},
r'_2(x) \pmod{x^2},
\cdots,
r'_{k+1}(x) \pmod{x^{k+1}}
\}
\]
covers all polynomials of degree $k$ or less, except for $r_0(x)$. Through inductive definition we obtain a set of congruence classes
\[
\{
r'_1(x) \pmod{x},
r'_2(x) \pmod{x^2},
\cdots,
r'_{n}(x) \pmod{x^{n}}
\}
\]
that covers all polynomials of degree less than $n$, except for $r_0(x)$.
\end{proof}
\section{Two Lemmas}~\label{lemmassec}
In this section, we prove two lemmas we'll need in the proof of Proposition~\ref{keyprop}. 
\begin{lemma}\label{oppo}
Suppose Theorem~\ref{2thm} is not true. Then for some $n$ there exist congruences $f(x) \equiv r_i(x) \pmod{b_i(x)}$, $i =1, 2, \cdots, n$, such that all the following conditions hold:
\begin{enumerate}[(1)]
    \item The set of congruences covers all polynomials of degree $n-1$ or less but not all polynomials.
    \item All the $b_i(x) = p_i(x)$'s are irreducible.
    \item If $k$ of the congruences have the same modulus $p(x)$, then $k < d$, where $d=\deg(p(x))$.
\end{enumerate}
\end{lemma}
\begin{proof}
Let $n$ be the smallest integer for which Theorem~\ref{2thm} does not hold. Hence there exist congruence classes $\langle f_1(x) \rangle + r_1(x)$, $\langle f_2(x) \rangle + r_2(x)$,  $\cdots$, and $\langle f_n(x) \rangle + r_n(x)$ that cover every polynomial of degree less than $n$, but not all polynomials in $\mathbb{F}_2[x]$. Let $\alpha(x)$ be a polynomial omitted by these congruence classes. We first show that we can adjust the congruence classes so that each $f_i(x)= p_i^{\lambda_i}(x)$ has only one irreducible factor. Suppose $f_i(x)=\prod\limits_{j=1}^{J} p_j(x)^{\lambda_j}$ with all $p_j(x)$'s irreducible. Then the following congruence equations
\[
\begin{cases}
\alpha(x) \equiv r_i(x) \pmod{p_1(x)^{\lambda_1}},
\\
\alpha(x) \equiv r_i(x) \pmod{p_2(x)^{\lambda_2}},
\\
\cdots,
\\
\alpha(x) \equiv r_i(x) \pmod{p_J(x)^{\lambda_J}}
\end{cases}
\]
cannot hold simultaneously. Therefore, there exists some $j$ such that $\alpha(x) \not\equiv r_i(x) \pmod{p_j(x)^{\lambda_j}}$. Now replace $f_i(x)$ with $p_j(x)^{\lambda_j}$, and we know the new set of congruence classes does not cover $\alpha(x)$ but still covers all polynomials of degree $n-1$ or less.

Then we move on to (3).  Fix an irreducible polynomial $p(x)$ of degree $d$, and let $\alpha_0(x)$ be the residue of $\alpha(x)$ modulo $p(x)$. There are 3 types of congruence classes:  
\[
\begin{cases}
r_i(x) \pmod{p(x)} \text{ such that }
\alpha(x) \not\equiv r_i(x) \pmod{p(x)}, \,\,\,\,\,\,\,\, 
(a), 
\\
\kappa(x) p(x) +\alpha_0(x) \pmod{p(x)^{\lambda}}
\text{ for some }
\lambda>1, 
\,\,\,\,\,\,\,\, 
(b), 
\\
r'_i(x) \pmod{q(x)} 
\text{ such that }
\gcd(p(x),q(x))=1, 
\,\,\,\,\,\,\,\, 
(c). 
\end{cases}
\]
For the last two types of congruence classes, they cover all polynomials $a(x)\equiv \alpha(x) \pmod{p(x)}$ with $\deg(a(x)) \leq n - 1$, but not $\alpha(x)$. So we rewrite them as
\[
\kappa_i p(x) + \alpha_0(x) \pmod{b_i(x) p(x)},
\]
where $b_i(x) = p(x)^{\lambda-1}$ for some $\lambda$ or $q(x)$ for some $(q(x),p(x))=1$, $1\leq i\leq n-k$. Note that $\deg(\alpha_0(x))< \deg(p(x))$, and we have a bijection between polynomials of degree $n-1-d$ or less and polynomials of degree $n-1$ or less in the form $\kappa_i(x) p(x) +\alpha_0(x)$. So the $n-k$ congruence classes $\kappa_i(x) \pmod{b_i(x)}$ cover all polynomials of degree $n-d-1$ or less, except for $\alpha(x)-\alpha_0(x)$. Recall Proposition~\ref{keyprop} holds for $n-k$, so we must have $n-k>n-d$, this is (3).

Next we discuss (2). By (3), there must be some $g(x) \equiv r_0(x) \pmod{p(x)}$ not included in any type-(a) congruence and $r_0(x) \neq \alpha_0(x)$. Let 
\[
M(x) = \prod \limits_{(q(x),p(x))=1 \text{ in a type-(c) congruence}} q(x),
\]
so $(M(x),p(x))=1$. Since type-(c) congruence classes omit some polynomials, they omit at least one congruence class $[g_0(x)]$. (If $g_0(x)$ is omitted, it is obvious $g_0(x)+M(x)$ is also omitted by all type-(c) congruences.) By Chinese remainder theorem, there exists a polynomial $g'(x)$ satisfying 
\[
\begin{cases}
g'(x) \equiv r_0(x) \pmod{p(x)}, 
\\
g'(x) \equiv g_0(x) \pmod{M(x)}. 
\end{cases}
\]
If we replace all type-(b) congruence classes with $\alpha_0(x) \pmod{p(x)}$, the new congruence classes still cover all polynomials of degree $n-1$ or less, although they cover $\alpha(x)$ now, they still do not cover $g'(x)$. So (1) still holds. By discussing the (3) property again, we can see (3) still holds. By repeating similar replacements for other irreducible factors we can assume (2). This completes the proof of the lemma. 
\end{proof} 

\begin{lemma}\label{inex2}
Let $S_1,S_2,\cdots,S_t$ be sets of polynomials in $\mathbb{F}_2[x]$ and $S_i$ consists of exactly $k_i$ congruence classes modulo $p_i(x)$, where $\deg(p_i)=d_i$ for $1 \leq i \leq t$ and $\gcd(p_i(x),p_j(x))=1$ for $1 \leq i\neq j \leq t$. Let $n$ be a positive integer and let $N$ be the number of polynomials of degree $n-1$ or less not included in any of the $S_i$'s. Let $s$ be any integer $1\leq s\leq t$. We have
\[
N>1+2^n\left(1-\sum\limits_{i=1}^s \frac{k_i}{2^{d_i}}\right)\prod\limits_{i=s+1}^t\left(1-\frac{k_i}{2^{d_i}}\right)-\left(1+\sum\limits_{i=1}^s k_i\right)\prod\limits_{i=s+1}^t(1+k_i).
\]
\end{lemma}
\begin{proof}  For any set $S$, let $C(S)$ be its characteristic function. Then we note that 
\[
1-\sum\limits_{i=1}^s C(S)\leq \prod\limits_{i=1}^s(1- C(S)).
\]
So the characteristic function of the set of polynomials not in any $S$ is:
\begin{align*}
&C
\left(
\overline {\bigcup\limits_{i=1}^t S_i}
\right) \\
&=
C
\left(
\bigcap\limits_{i=1}^t \overline S_i
\right) \\
&=
\prod\limits_{i=1}^t C(\overline S_i) \\
&=
\prod\limits_{i=1}^t (1-C(S_i)) \\
&=
\prod\limits_{i=1}^s (1-C(S_i)) \prod\limits_{i=s+1}^t (1-C(S_i)) \\
&\geq(1-\sum\limits_{i=1}^s C(S_i)) \prod\limits_{i=s+1}^t (1-C(S_i)) \\
&=
1-\sum_{i=1}^t C(S_i)+\sum\limits_{i,j}{}'C(S_i\cap S_j)-\cdots, 
\end{align*}
where $\sum'$ indicates that at most one subscript is less than or equal to $s$.

If $d_i\leq n$, then there are $k_i$ elements of $S_i$ among polynomials of degrees at most $d_i-1$, and $2^{d-d_i} \cdot k_i$ more of degree $d \geq d_i$. So we have 
\[
\sum\limits_{\deg(r(x))<n}C(S_i)(r(x))=k_i
\left(
1+\sum\limits_{d=d_i}^{n-1}2^{d-d_i}
\right)
=2^{n-d_i} \cdot k_i.
\]
If $d_i>n$, then there are at least 0 and at most $k_i$ elements of $S_i$ among polynomials whose degrees are less than $n$, or we have 
\[
0\cdot k_i\leq\sum\limits_{\deg(r)<n}C(S_i)(r(x))<0\cdot k_i+k_i.
\]
To put the cases together, we have
\[
\sum\limits_{deg(r)<n}C(S_i)(r)=2^{n-d_i}k_i+E_i,
\]
where $|E_i|<k_i$. By the Chinese remainder theorem, there are at most $k_i k_j \cdots k_z$ polynomials in $S_i \cap S_j \cap \cdots \cap S_z$ of degree less than $d_i+d_j+ \cdots d_z$, so 
\[
\sum\limits_{\deg(r)<n}C(S_i\cap S_j\cap \cdots \cap S_z)(r(x))=2^{n-d_i-d_j- \cdots - d_z}k_ik_j \cdots k_z+E_{ij \cdots z}, 
\]
where $|E_{ij \cdots z}|<k_ik_j \cdots k_z$. Then 
\begin{align*}
N
&=
\sum\limits_{\deg(r(x))<n}C
\left(
\overline {\bigcup\limits_{i=1}^t S_i}
\right)
(r(x)) \\
&\geq \sum\limits_{\deg(r)<n}
\left(
1-\sum_{i=1}^t C(S_i)+\sum\limits_{i,j}{}'C(S_i\cap S_j)-\cdots
\right)
(r(x)) \\
&=
2^n-\sum_{i=1}^t2^{n-d_i}k_i+ \sum\limits_{i,j}{}'2^{n-d_i-d_j}k_ik_j-\cdots+E \\
&=
2^n
\left(
1-\sum\limits_{i=1}^s \frac{k_i}{2^{d_i}}
\right)
\prod\limits_{i=s+1}^t
\left(
1-\frac{k_i}{2^{d_i}}
\right)+E, 
\end{align*}
where 
\[
|E|<
\left|
\sum_{i=1}^tk_i- \sum\limits_{i,j}{}'k_ik_j+\cdots
\right|
=
\left(1+\sum\limits_{i=1}^s k_i\right)
\prod\limits_{i=s+1}^t(1+k_i), 
\]
which completes the proof of the lemma.
\end{proof}
\section{Proof of Proposition~\ref{keyprop}}\label{proofsec}
In this section, we continue to prove Proposition~\ref{keyprop}.  By Lemma \ref{inex2} it suffices to prove for any $n\geq 10$ and some $s$ we have 
\begin{equation}\label{ineq}
2^n
\left(
1-\sum\limits_{i=s+1}^t \frac{k_i}{2^{d_i}}
\right)
\prod\limits_{i=1}^s
\left(
1-\frac{k_i}{2^{d_i}}
\right) 
\geq 
\left(
1+\sum\limits_{i=s+1}^tk_i
\right)
\prod\limits_{i=1}^s(1+k_i). 
\end{equation}, where each set $S_i$ here is the union of all congruence classes modulo $p_i(x)$ that appeared in set of congruences we found through Lemma \ref{oppo}.

Without loss of generality, we can assume $1<d_1 \leq d_2 \leq \cdots\leq d_t$.
Note that the right side of Inequality~\eqref{ineq} is the product of $(s+1)$ factors with sum $n+s+1$. Let $s=\min(\lceil{\frac{n}{3}-1}\rceil,t-1)$. Then the right side of Inequality~\eqref{ineq} is at most
\[
\left(
\frac{n+s+1}{s+1}
\right)^{s+1}
\leq 
\left(
\frac{n+\frac{n}{3}}{\frac{n}{3}}
\right)^{\frac{n}{3}}
=2^{\frac{2n}{3}}.
\]
Now we work on the left side of Inequality~\eqref{ineq}, which we denote by $L$. We want to prove
\[ L \geq 2^{\frac{2n}{3}}. \].
To handle the first factor, we need the following lemma. 
\begin{lemma}\label{fubini}
Let d be a positive integer. The number of irreducible polynomials in $\mathbb{F}_2[x]$ whose degree is greater than 1 and at most d, denoted as $N^2_d$, satisfies the inequality $N^2_d \leq 3 \cdot 2^{d-3}$.
\end{lemma}
\begin{proof} Note that for $d\in \mathbb{N}$ there are $n^2_d=\frac{1}{d}\sum\limits_{k|d}\mu(k)2^{\frac{d}{k}}$ $d$-th degree irreducible polynomials in $\mathbb{F}_2[x]$. Hence $n^2_d\leq \frac{2^d}{d}\leq\frac{3\cdot2^d}{16}$ for $d>5$. For $1\leq d\leq 5$ one can directly check that Lemma~\ref{fubini} holds, so $N^2_d=N^2_5+\sum\limits_{i=6}^dn^2_i \leq N^2_5+\sum\limits_{i=6}^d\frac{3\cdot2^i}{16}\leq \sum\limits_{i=0}^d\frac{3\cdot2^i}{16}=\frac{3\cdot2^d}{8}$ for $n>5$, thus the lemma is true for all $d\in \mathbb{N}$.
\end{proof}
By Lemma \ref{fubini}, for any $n\in \mathbb{N}$, we know $n \leq N^2_{d_n}\leq \frac{3\cdot2^{d_n}}{8}$ (since $p_i$'s here are irreducible, their number cannot exceed the number of all irreducible polynomials whose degrees are no more than $d_n$), so we have $d_n\geq \log_2(\frac{8n}{3})$ and $2^{d_n}\geq \frac{8n}{3}$.
Therefore, we have 
\[
1-\sum\limits_{i=s+1}^t \frac{k_i}{2^{d_i}}\geq 1-\sum\limits_{i=[\frac{n}{3}]+1}^t \frac{k_i}{2^{d_[\frac{n}{3}]}}\geq1-\frac{\frac{2n}{3}}{\frac{8}{3}\cdot\frac{n}{3}}=\frac{1}{4}.
\]
We now apply the adjustment method to handle the remaining factors.
First, if there  exist $p_1(x), p_2(x)$ such that $d_1=d_2$ and $1\leq k_1\leq k_2<d_2-1$, we know
\[
\left(
1-\frac{k_1}{2^{d_1}}
\right)
\left(
1-\frac{k_2}{2^{d_2}}
\right)
-
\left(
1-\frac{k_1-1}{2^{d_1}}
\right)
\left(
1-\frac{k_2+1}{2^{d_2}}
\right)
=
\frac{k_2-k_1+1}{2^{2d_1}}>0,
\]
and hence
\[
\left(
1-\frac{k_1}{2^{d_1}}
\right)
\left(
1-\frac{k_2}{2^{d_2}}
\right)
>
\left(
1-\frac{k_1-1}{2^{d_1}}
\right)
\left(
1-\frac{k_2+1}{2^{d_2}}
\right).
\]
On the other hand, if there  exist $p_1(x), p_2(x)$ such that $d_1<d_2$ and $0\leq k_1<d_1-1$, then
\begin{align*}
&
\left(
1-\frac{k_1}{2^{d_1}}
\right)
\left(
1-\frac{k_2}{2^{d_2}}
\right)
-(
\left(
1-\frac{k_1+1}{2^{d_1}}
\right)
\left(
1-\frac{k_2-1}{2^{d_2}}
\right) \\
&=\frac{2^{d_2}-2^{d_1}+k_1-k_2+1}{2^{d_1+d_2}} \\
&>\frac{2^{d_2-1}-k_2}{2^{d_1+d_2}} \\
&>\frac{2^{d_2-1}-d_2}{2^{d_1+d_2}} \\
&> 0. 
\end{align*}
Therefore, we have 
\[
\left(
1-\frac{k_1}{2^{d_1}}
\right)
\left(
1-\frac{k_2}{2^{d_2}}
\right)
>
\left(
1-\frac{k_1+1}{2^{d_1}}
\right)
\left(
1-\frac{k_2-1}{2^{d_2}}
\right). 
\]
By the discussion above, we can assume that the remaining factors ``fill" the irreducible moduli in an increasing order of degrees, from the lowest degree to the highest. Since $k_i<d_i$, we have 
\begin{align*}
&\prod\limits_{i=1}^s
\left(
1-\frac{k_i}{2^{d_i}}
\right) \\
&\geq
\left(
1-\frac{1}{2^2}
\right)
\left(
1-\frac{2}{2^3}
\right)
\left(
1-\frac{2}{2^3}
\right)^{[\frac{[\frac{n}{3}]-1-3}{2}]+1} \\
&>
\left(
1-\frac{1}{2^2}
\right)
\left(
1-\frac{2}{2^3}
\right)
\left(1-\frac{2}{2^3}
\right)^{\frac{\frac{n}{3}-4}{2}+1}. 
\end{align*}
So for $n\geq 10$, we have 
\begin{align*}
&2^n
\left(
1-\sum\limits_{i=s+1}^t \frac{k_i}{2^{d_i}}
\right)
\prod\limits_{i=1}^s
\left(
1-\frac{k_i}{2^{d_i}}
\right) \\
&\geq 
\left(1-\frac{1}{4}\right)
\left(1-\frac{1}{2^2}\right)
\left(1-\frac{2}{2^3}\right)
\left(1-\frac{2}{2^3}\right)^{\frac{\frac{n}{3}-4}{2}+1}\\
&=
2^{\frac{2n}{3}-4}\cdot3^{\frac{n}{6}+1} \\
&>
2^{\frac{2n}{3}} \\
&\geq \left(
1+\sum\limits_{i=s+1}^t k_i \right) \prod\limits_{i=1}^s(1+k_i).
\end{align*}
This completes the proof of Proposition~\ref{keyprop}.

\section{Proof of Theorem 1}~\label{futureresultssec}
It remains to prove Theorem~\ref{2thm} for $n<10$. In this section our set of congruence classes discussed is always in the form mentioned in Lemma~\ref{oppo}, since this lemmas always hold. We state some facts to help us handle these cases. On the one hand, recall that there are at most $\lceil 2^{n-d_i}\rceil k_i$ polynomials whose degree is at most $n-1$ in $S_i$ and they must cover \{$f(x) \in \mathbb{F}_2[x]: \deg(f(x))\leq n-1$\}. So for a set of congruence classes satisfying Lemma \ref{oppo}, when $n=1$, there are at most 1 constant covered; for $n=2$ at most 2 polynomials of degree at most 1 covered; and for $n=3$ at most $\frac{2^3}{2^2}+2=4$ polynomials covered. If $4\leq n\leq 8$, there are at most $\left(\frac{1}{4}+\frac{2}{8}+\frac{2}{8}+\frac{3}{16}\right)2^n=\frac{15}{16}\cdot2^n$ polynomials covered, a contradiction.

For $n=9$, there are at most $\left(\frac{1}{4}+\frac{2}{8}+\frac{2}{8}+\frac{3}{16}+\frac{1}{16}\right)2^n=2^n$ polynomials covered. But to make the equality hold, the $S_i$'s must be disjoint. But by the Chinese remainder theorem this is impossible, a contradiction. This completes the proof of Theorem~\ref{2thm}.
\section{Proof of Theoerm~\ref{qthm} and some discussion}~\label{generalizedresultssec}
To prove Theorem~\ref{qthm} we need to prove a similar lemma in $\mathbb{F}_q[x]$. 
\begin{lemma}\label{inexq}
Let $S_1,S_2,\cdots,S_t$ be sets of polynomials in $\mathbb{F}_q[x]$ and $S_i$ consist of exactly $k_i$ congruence classes mod $p_i(x)$, where $\deg(p_i(x))=d_i$ for $1 \leq i \leq t$ and $\gcd(p_i(x),p_j(x))=1$ for $1 \leq i \neq j \leq t$. Let $n$ be a positive integer and $N$ be the number of polynomials of degree $n-C-1$ or less and not included in any of the $S_i$'s, and let $s$ be any integer $1\leq s\leq t$, here $C$ is a constant. Then,
\[
N
>
1+q^{n-C}
\left(
1-\sum\limits_{i=1}^s \frac{k_i}{q^{d_i}}
\right)
\prod\limits_{i=s+1}^t
\left(
1-\frac{k_i}{q^{d_i}}
\right)
-
\left(1+\sum\limits_{i=1}^s k_i
\right)
\prod\limits_{i=s+1}^t(1+k_i).
\]
\end{lemma}
The proof is similar to that of Lemma \ref{inex2} using the inclusion-exclusion principle. 

Now we prove Theorem~\ref{qthm}. Note that in the proof of Lemma \ref{oppo} we did not use any property specific in $\mathbb{F}_2[x]$ from general $\mathbb{F}_q[x]$, and we can still make the discussion. So, for $n\geq 10$, the inequality we need to prove Theorem 2 is now
\[
q^n(1-\sum\limits_{i=s+1}^t \frac{k_i}{q^{d_i}})\prod\limits_{i=1}^s(1-\frac{k_i}{q^{d_i}})\geq 4^\frac{n}{3}.
\]
The adjustment method can still be applied. Note that $N^q_d\leq \frac{q-1}{2q}q^{d}$ and $q^{d_n}>\frac{2qn}{q-1}$, where $N^q_d$ is the number of irreducible polynomials in $\mathbb{F}_q[x]$ whose degree is equal or less than $d$ but no less than 2. Again, the proof is similar to that of Lemma~\ref{fubini}. Hence
\[
q^n
\left(
1-\sum\limits_{i=s+1}^t \frac{k_i}{q^{d_i}}
\right)
\prod\limits_{i=1}^s
\left(
1-\frac{k_i}{q^{d_i}}
\right)
\geq q^{n-1}\prod\limits_{i=1}^s\frac{8}{9}.
\]
For $2\leq n\leq 9$, since
\[
\frac{\lceil q^{n-d_i}\rceil k_i}{q^n}\leq \frac{1}{9}, 
\]
we know Theorem~\ref{qthm} follows. 

Next, we continue to discuss Theorems~\ref{2thm} and \ref{qthm} and propose a conjecture. 

First, by appropriately replacing Lemmas~\ref{oppo} and~\ref{inex2}, we can prove that every set of $n$ congruence equations covering all polynomials of degree exact $n$ (or any $m>n$) covers the polynomial ring.

Second, according to our results, when $q=3$ and $n=2$, if a set of 2 congruence classes in $\mathbb{F}_3[x]$ covers all polynomials whose degrees are at most 1, they cover the polynomial ring. However, one may note that if a set of 2 congruence classes covers all polynomials of degree at most 0, it already covers the whole ring. The proof of this is similar to how we handle the case $q=2$ and $n\leq 9$. By testing on small $q$'s and $n$'s, we conjecture the lower bound of degree can be $n/(q-1)$. This leads us to the conjecture in Section~\ref{introsec}. 
\begin{conjecture}
Let $f_1(x), f_2(x), \cdots, f_n(x), r_1(x), r_2(x), \cdots, r_n(x)$ be polynomials in $\mathbb{F}_q[x]$ satisfying $\deg(r_i(x)) < \deg(f_i(x))$ for $1\leq i\leq n$. Suppose that there exists some polynomial $g(x) \in \mathbb{F}_q[x]$ that satisfies none of the congruences
\[
g(x) \equiv r_i(x) \pmod{f_i(x)}, 
\]
for all $1 \leq i \leq n$. 
Then there exists a polynomial $g_0(x) \in \mathbb{F}_q[x]$ that satisfies none of the congruence equations above and $\deg(g_0(x)) < {n}/{(q - 1)}$. 
\end{conjecture}

Though, even though the conjecture is very likely to be true, it is unlikely that this method we used above can be used to prove the desired results. The reason is that, to apply a similar method, we need to restrict the asymptotic approximation of 
$q^{n/(q-1)} (1-\sum\limits_{i=s+1}^t \frac{k_i}{q^{d_i}}) \prod\limits_{i=1}^{s}(1-\frac{k_i}{q^{d_i}})$ to at least $2^{2n/3}$, which seems to be impossible since the first factor $q^{n/{q-1}}$ is already too small.

\bibliographystyle{elsarticle-num} 
\bibliography{bib.bib}



\end{document}